\documentclass[12pt]{amsart}
\usepackage{graphicx}
\usepackage{amsmath}
\usepackage{amssymb}
\usepackage{latexsym}
\usepackage{amsfonts}
\usepackage{a4wide}
\usepackage{longtable}

\newcommand{\Z}{\mathbb{Z}}
\newcommand{\Q}{\mathbb{Q}}

\renewcommand{\P}{\mathbb{P}}

\newcommand{\lra}{\longrightarrow}

\title{Some monodromy groups\\ of finite index in $Sp_4(\Z)$}
\author{J\"org Hofmann, Duco van Straten}

\begin{document}
\begin{abstract}
We determine the index of five of the $7$  hypergeometric
Calabi-Yau operators that have finite index in $Sp_4(\Z)$ and
in two cases give a complete description of the monodromy group. 
Furthermore we found six non-hypergeometric Calabi-Yau operators 
with finite index in $Sp_4(\Z)$, most notably a case where the 
index is one.
\end{abstract}

\maketitle

\section{Introduction}
The fourteen hypergeometric fourth order operators related to mirror
symmetry for complete intersections in weighted projective space have
always been treated as a single group, with very similar properties. An
explicit description of monodromy matrices has been known since a long 
time. It came therefore as a surprise to us that recently S. Singh and 
T. N. Venkataramana showed that in at least three of the fourteen cases the
monodromy is of finite index in $Sp_4(\Z)$. On the other hand, the work 
of C. Brav and  H. Thomas showed that in at least $7$ of the $14$ cases 
the monodromy is of infinite index. In a further paper, S. Singh has shown 
that the monodromy is finite in the four remaining cases. 
So an interesting dichotomy has arisen in the class of Calabi-Yau operators. 
In this note we give a precise determination of two of the groups of finite 
index and determine the index in three more cases. Furthermore, six 
non-hypergeometric Calabi-Yau operators are identified which have finite 
index in $Sp_4(\Z)$. 

\section{The fourteen hypergeometric families}

The general quintic hypersurface in ${\bf P}^4$ and the remarkable
enumerative properties of the Picard-Fuchs operator of the mirror family 
\[\theta^4-5^5 x (\theta+\frac{1}{5})(\theta+\frac{2}{5})(\theta+\frac{3}{5})(\theta+\frac{4}{5})\]
discovered by {\sc Candelas, de la Ossa, Green} and {\sc Parkes} \cite{COGP} stands at the beginning of much of the interest in the mirror symmetry phenomenon that continues up to the present day.
The above example was readily generalised to the case of smooth Calabi-Yau threefolds
in {\em weighted} projective space, producing three further cases, \cite{M}, \cite{KT1}. Then
{\sc Libgober and Teitelbaum} \cite{LT} produced mirror families for the other four Calabi-Yau 
complete intersections in ordinary projective spaces. A final generalisation consisted of looking at smooth complete intersections Calabi-Yau threefolds in weighted projective spaces, leading to a further five cases, \cite{KT2}. In all these $13$ cases the Picard-Fuchs operator is 
{\em hypergeometric} and takes the form
\[\theta^4-Nz(\theta+\alpha_1)(\theta+\alpha_2)(\theta+\alpha_3)(\theta+\alpha_4)\] 
It was remarked by several authors that in fact there is an {\em overlooked, fourteenth} case,
corresponding to the complete intersection of hypersurfaces of degree $2$ and $12$ in
 $\P^5(1,1,1,1,4,6)$, which represents a Calabi-Yau threefold with a singularity, \cite{A}, \cite{DM}, \cite{R}. 
From the point of view of differential equations the fourteen hypergeometric equations are 
characterised as fourth order hypergeometrics with exponents $0,0,0,0$ at $0$ that carry a {\em monodromy invariant lattice}. This leads to a monodromy group that is (conjugate to) a sub-group of 
$Sp_4(\Z)$ and a necessary (and, after the fact sufficient) condition for this to happen is 
that the characteristic  polynomial of the monodromy around $\infty$ is a product of cyclotomic polynomials, which leads immediately to the $14$ cases. Below we summarise the situation in a table. The last column give the number as it appears in the table \cite{AESZ}.

{\small
\[
\begin{array}{|c|c|c|c|}
\hline
\textup{Case}&N&\alpha_1,\alpha_2,\alpha_3,\alpha_4&AESZ\\[1mm]
\hline
\hline
&&&\\[-2mm]
\P^4[5]&5^5& \frac{1}{5},\frac{2}{5}, \frac{3}{5},\frac{4}{5}&1\\[2mm]
\hline
&&&\\[-2mm]
\P^4(1,1,1,1,2)[6]&2^43^6&\frac{1}{6},\frac{1}{3}, \frac{2}{3},\frac{5}{6}&8\\[2mm]
\P^4(1,1,1,1,4)[8]&2^{16}&\frac{1}{8},\frac{3}{8}, \frac{5}{8},\frac{7}{8}&7\\[2mm]
\P^4(1,1,1,2,5)[10]&2^85^5&\frac{1}{10},\frac{3}{10}, \frac{7}{10},\frac{9}{10}&2\\[2mm]
\hline
&&&\\[-2mm]
\P^5[3,3]&3^6&\frac{1}{3},\frac{1}{3}, \frac{2}{3},\frac{2}{3}&4\\[2mm]
\P^5[2,4]&2^{10}&\frac{1}{4},\frac{2}{4}, \frac{3}{4},\frac{3}{4}&6\\[2mm]
\P^6[2,2,3]&2^43^3&\frac{1}{3},\frac{1}{2}, \frac{1}{2},\frac{2}{3}&5\\[2mm]
\P^7[2,2,2,2]&2^8&\frac{1}{2},\frac{1}{2}, \frac{1}{2},\frac{1}{2}&3\\[2mm]
\hline
&&&\\[-2mm]
\P^5(1,1,1,1,2,2)[4,4]&2^{12}&\frac{1}{4},\frac{1}{4}, \frac{3}{4},\frac{3}{4}&10\\[2mm]
\P^5(1,1,1,1,1,2)[3,4]&2^63^3&\frac{1}{4},\frac{1}{3}, \frac{2}{3},\frac{3}{4}&11\\[2mm]
\P^5(1,1,1,2,2,3)[4,6]&2^{10}3^3&\frac{1}{6},\frac{1}{4}, \frac{3}{4},\frac{5}{6}&12\\[2mm]
\P^5(1,1,2,2,3,3)[6,6]&2^{8}3^6&\frac{1}{6},\frac{1}{6}, \frac{5}{6},\frac{5}{6}&13\\[2mm]
\P^5(1,1,1,1,1,3)[2,6]&2^83^3&\frac{1}{6},\frac{1}{2}, \frac{1}{2},\frac{5}{6}&14\\[2mm]
\hline
&&&\\[-2mm]
\P^5(1,1,1,1,4,6)[2,12]&2^{12}3^6&\frac{1}{12},\frac{5}{12},\frac{7}{12},\frac{11}{12}&9\\[2mm]
\hline
\end{array}
\]
}
The factor $N$ is introduced to make the power series expansion around $0$ of the
holomorphic solution have integral coefficients in a minimal way. We call $N$ the {\em discriminant} of the operator; the critical point is then located at $x=1/N=:x_c$. In terms of the exponents $\alpha_1,\alpha_2,\alpha_3=1-\alpha_2,
\alpha_4=1-\alpha_1$ it can be given as (see \cite{B})
\[N= \prod_{i=1}^4 N(\alpha_i)\]
where
\[N(\frac{r}{s}):=m(s);\;\;\; m(s):=s\prod_{p \mid s} s^{1/p-1}\]
so
\[
\begin{array}{|c|c|c|c|c|c|c|c|c|}
\hline
s&2&3&4&5&6&8&10&12\\
\hline 
&&&&&&&&\\[-2mm]
m(s)&2^2&3^{3/2}&2^3&5^{5/4}&2^23^{3/2}&2^4&2^25^{5/4}&2^33^{3/2}\\
\hline
\end{array}
\]

\section{Monodromy matrices}
The explicit description of the monodromy of the general hypergeometric operator \[(\theta+\beta_1-1)\ldots (\theta+\beta_n-1)-x(\theta+\alpha_1) \ldots (\theta+\alpha_n)\]
has a long history. In his thesis, \cite{L} {\sc Levelt} showed the existence of a basis where the monodromy 
around $\infty$ and $0$ are given by the companion matrices of the characteristic polynomials 
\[f(T)=\prod_{k=1}^n(T-e^{2\pi i\alpha_k}),\;\;g(T)=\prod_{k=1}^n(T-e^{2\pi i\beta_k})\]
However, for our purpose it is natural to work with other bases.
First of all, for all our operators there is a unique {\em Frobenius basis} of
solutions around $0$ of the form 
\[
\begin{array}{rcl}
\Phi_0(x)&=&f_0(x)\\[1mm]
\Phi_1(x)&=&\log(x)f_0(x)+f_1(x)\\[1mm]
\Phi_2(x)&=&\frac{1}{2}\log(x)^2f_0(x)+\log(x)f_1+f_2(x)\\[1mm]
\Phi_3(x)&=&\frac{1}{6}\log(x)^3f_0(x)+\frac{1}{2}\log^2(x)f_2(x)+\log(x)f_1(x)+f_3(x)\\[1mm]
\end{array}
\]
where $f_0=1+\ldots \Z[[x]]$ and $f_1,f_2,f_3 \in x\Q[[x]]$.
The basis of solutions
\[ y_k(x):=\frac{1}{(2\pi i)^k}\Phi_k(x)\]
is called the {\em normalised Frobenius basis}; the monodromy around 
$0$ in this basis is given by
\[
M_F=\left(
\begin{array}{cccc}
1&1&\frac{1}{2}&\frac{1}{6}\\
0&1&1&\frac{1}{2}\\
0&0&1&1\\
0&0&0&1\\
\end{array}
\right)
\]
In this basis the monodromy invariant symplectic form is given by
\[
S_F=\left(
\begin{array}{cccc}
0&0&0&1\\
0&0&-1&0\\
0&1&0&0\\
-1&0&0&0\\
\end{array}
\right)
\]
and the monodromy around $x_c$ is a symplectic reflection
\[ v \lra v - \frac{1}{d}\langle C,v \rangle C\]
in a vector $C$ that represents the vanishing cycle and
which has the form 
\[C=(d,0,b,a)\]
where $d:=H^3$ is the {\em degree} of the ample generator,
$b:=c_2(X)H/24$ and $a: =\lambda c_3(X)$ are the characteristic numbers of the 
corresponding Calabi-Yau threefold $X$ and
\[ \lambda=\frac{\zeta(3)}{(2\pi i)^3}\,.\]
A further important invariant is the number
\[k= \frac{H^3}{6}+\frac{c_2(X)\cdot H}{12}=\frac{d}{6}+2b\]
which is equal to the dimension of the linear system $|H|$.

The base-change by the matrix
\[A= \left(\begin{array}{cccc}
0&0&1&0\\
0&0&0&1\\
0&d&d/2&-b\\
-d&0&-b&-a\\
\end{array} \right)\]
conjugates the matrices $M_F$ and $N_F$ to 
\[
M:=AM_FA^{-1} =\left( \begin{array}{cccc}
 1&1&0&0\\
 0&1&0&0\\
d&d&1&0\\
0&-k&-1&1\\
\end{array}\right),\;\;
N:=A N_FA^{-1}=\left( \begin{array}{cccc}
 1&0&0&0\\
 0&1&0&1\\
0&0&1&0\\
0&0&0&1\\
\end{array}\right)
\]
which are now in the  integral symplectic group 
\[Sp_4(\Z)=\{M\;|\;M^t \cdot S \cdot M=I \}\]
realisd as set of integral matrices that preserve the {\em standard symplectic
form}
\[S:=\left( \begin{array}{cccc}
 0&0&1&0\\
 0&0&0&1\\
-1&0&0&0\\
0&-1&0&0\\
\end{array}\right)\,.\]
This is the form of the generators that can be found in \cite{CYY}.

So the monodromy group $G(d,k)$ of the differential operator is
the group generated by these two matrices $M$ and $N$.
It was observed in \cite{CYY} that the monodromy group in fact
is contained in a congruence subgroup
\[ G(d,k) \subset \Gamma(d,gcd(d,k))\]
where
$\Gamma(d_1,d_2)$, $d_2 \mid d_1$,  consist of those matrices $A$ in $Sp_4(\Z)$ for which
\[A \equiv \left( \begin{array}{cccc}1&*&*&*\\0&*&*&*\\0&0&1&0\\0&*&*&*\end{array} \right) \mod d_1,\;\;\; A \equiv \left( \begin{array}{cccc}1&*&*&*\\0&1&*&*\\0&0&1&0\\0&0&*&1 \end{array} \right) \mod d_2  \]
The index of this group in $Sp_4(\Z)$ was computed by {\sc C. Erdenberger} \cite{CYY}, Appendix,
as
\[|Sp_4(\Z):\Gamma(d_1,d_2)|=d_1^4\prod_{p \mid d_1}(1-p^{-4})d_2^2\prod_{p \mid d_2}(1-p^{-2}), \]
where the product runs over the primes dividing $d_1$ resp. $d_2$.

The parameters $(d,k)$ suggest a natural way to order the list of
$14$ hypergeometric cases. Remarkably, this ordering coincides with the one
obtained by either using the first instanton number $n_1$ (rational curves of degree one) or the discriminant $N$.
\[
\begin{array}{|c|c|c|c|c|r|r|c|}
\hline
(d,k)&\alpha_1,\alpha_2&H^3&c_2\cdot H & c_3&n_1&N&AESZ\\
\hline
\hline
(1,4)&\frac{1}{12},\frac{5}{12}&1&46&-484&678816&2985984&9\\[1mm]
\hline
(1,3)&\frac{1}{10},\frac{3}{10}&1&34&-288&231200&800000&2\\[1mm]
\hline
(1,2)&\frac{1}{6},\frac{1}{6}&1&22&-120&67104&86624&13\\[1mm]
\hline
(2,4)&\frac{1}{8},\frac{3}{8}&2&44&-296&29504&65536&7\\[1mm]
\hline
(2,3)&\frac{1}{6},\frac{1}{4}&2&32&-156&15552&27648&12\\[1mm]
\hline
(3,4)&\frac{1}{6},\frac{1}{3}&3&42&-204&7884&11664&8\\[1mm]
\hline
(4,5)&\frac{1}{6},\frac{1}{2}&4&52&-256&4992&6912&14\\[1mm]
\hline
(4,4)&\frac{1}{4},\frac{1}{4}&4&40&-144&3712&4096&10\\[1mm]
\hline
(5,5)&\frac{1}{5},\frac{2}{5}&5&50&-200&2875&3125&1\\[1mm]
\hline
(6,5)&\frac{1}{4},\frac{1}{3}&6&48&-156&1944&1728&11\\[1mm]
\hline
(8,6)&\frac{1}{4},\frac{1}{2}&8&56&-176&1280&1024&6\\[1mm]
\hline
(9,6)&\frac{1}{3},\frac{1}{3}&9&54&-144&1053&729&4\\[1mm]
\hline
(12,7)&\frac{1}{3},\frac{1}{2}&12&60&-144&720&432&5\\[1mm]
\hline
(16,8)&\frac{1}{2},\frac{1}{2}&16&64&-128&512&256&3\\[1mm]
\hline
\end{array}
\]
We remark further that the invariants $d$ and $k$ can be expressed
directly in terms of the defining exponents $\alpha_1,\alpha_2$ as follows:
\[d=4(1-\cos(2\pi \alpha_1))(1-\cos(2\pi \alpha_2)),k=4-2\cos(2\pi \alpha_1)-2\cos(2\pi \alpha_2)\]
which can be expressed as saying that 
\[2-2\cos(2\pi \alpha_1)\;\; \textup{and}\;\; 2- 2\cos(2\pi \alpha_2)\] 
are roots of the quadratic polynomial $X^2-kX+d=0$.

\section{Results}

During the last year important progress has been made in understanding 
the nature of the monodromy group $G(d,k)$.\\
 
{\bf Theorem 1} (C. Brav and H. Thomas, \cite{BT})

The group $G(k,d)$ has {\em infinite index} for the seven pairs
\[(d,k)=(1,4),(2,4),(4,5),(5,5),(8,6),(12,7),(16,8)\]

{\bf Theorem 2} (S. Singh and T. Venkataramana, S. Singh, \cite{SV}, \cite{S})

The group $G(k,d)$ has {\em finite index} for the other seven pairs
\[(d,k)=(1,3),(1,2),(2,3),(3,4),(4,4),(6,5),(9,6)\]

To these results we add\\

{\bf Theorem 3} The index $|Sp_4(\Z):G(d,k)|$ is given
by the following table
\[
\begin{array}{|c|c|c|c|c|c|c|c|}
\hline
(d,k)&(1,3)&(1,2)&(2,3)&(3,4)&(4,4)&(6,5)&(9,6)\\
\hline
&&&&&&&\\[-3mm]
\textup{Index $G(d,k)$ }&6&10&960&2^93^55^2&2^{20}3^25&2^{10} 3^6 5^2(?)&2^8 3^{13} 5^2 (?)\\
\hline
\textup{Index $\Gamma(d,gcd(d,k))$}&1&1&15&2^4 5&2^63^25&2^4 3^1 5^2&2^7 3^4 5\\
\hline
\end{array}
\]

The index of the last two entries is at least as big as the number indicated.
For easy comparison we also included the index of the corresponding group $\Gamma(d,gcd(d,k))$
in $Sp_4(\Z)$.\\
 
On the first two groups we can be very precise:\\

{\bf Theorem 4}

(i) The group $G(1,3)$ of index $6$ in $Sp_4(\Z)$
is exactly the group of matrices $A \in Sp_4(\Z)$ with the property 
the that $A \mod 2$
preserves the five-tuple of vectors of $(\Z/2)^4$
\[\{ \left( \begin{array}{c}0\\0\\0\\1 \end{array} \right),\;
\left( \begin{array}{c}0\\1\\0\\1 \end{array} \right),\;
\left( \begin{array}{c}0\\1\\1\\0 \end{array} \right),\;
\left( \begin{array}{c}1\\1\\0\\0 \end{array} \right),\;
\left( \begin{array}{c}1\\1\\1\\0 \end{array} \right)\}
\]

(ii) The group $G(1,2)$ of index $10$ in $Sp_4(\Z)$
is exactly the group of matrices $A \in Sp_4(\Z)$ with the property 
the that $A \mod 2$ preserves the pair of triples of vectors of $(\Z/2)^4$
\[\{ \{\left( \begin{array}{c}0\\0\\1\\0 \end{array} \right),\;
\left( \begin{array}{c}1\\1\\0\\0 \end{array} \right),\;
\left( \begin{array}{c}1\\1\\1\\0 \end{array} \right)\}
,\{
\left( \begin{array}{c}0\\0\\1\\1 \end{array} \right),\;
\left( \begin{array}{c}0\\1\\0\\0 \end{array} \right),\;
\left( \begin{array}{c}0\\1\\1\\1 \end{array} \right)
\}\}
\]

\section{Explanation of Theorem 3 and 4}

In order to determine the index of a sub-group in a given group, there is
the classical method of {\sc Todd} and {\sc Coxeter} called {\em coset-enumeration}. 
This has been developed into an effective computational tool that is 
implemented in {\tt GAP}, \cite{GAP}, the main tool for computational group 
theory. For details on this circle of ideas we refer to \cite{Neu}.

For this to work one needs a good presentation of $Sp_4(\Z)$ in terms of 
{\em generators} and {\em relations}. We used a presentation of
$Sp_4(\Z)$ described by {\sc Behr} in \cite{Be}, that uses $6$ generators and $18$ relations, and that is based on the root system for the symplectic group.
The six generating matrices are:
\begin{align*}
 x_{\alpha} & = \begin{pmatrix}
            1  & 1  & 0  & 0\\
            0  & 1  & 0  & 0\\
            0  & 0  & 1  & 0\\
            0  & 0  & -1 & 1\\
           \end{pmatrix},\;
           x_{\beta} = \begin{pmatrix}
            1  & 0  & 0  & 0\\
            0  & 1  & 0  & 1\\
            0  & 0  & 1  & 0\\
            0  & 0  & 0  & 1\\
           \end{pmatrix},\;x_{\alpha+\beta} = \begin{pmatrix}
            1  & 0  & 0  & 1\\
            0  & 1  & 1  & 0\\
            0  & 0  & 1  & 0\\
            0  & 0  & 0  & 1\\
           \end{pmatrix}\\
x_{2\alpha+\beta} &= \begin{pmatrix}
            1  & 0  & 1  & 0\\
            0  & 1  & 0  & 0\\
            0  & 0  & 1  & 0\\
            0  & 0  & 0  & 1\\
           \end{pmatrix},\;
 w_{\alpha} = \begin{pmatrix}
            0  & -1 & 0  & 0\\
            1  & 0  & 0  & 0\\
            0  & 0  & 0  & -1\\
            0  & 0  & 1  & 0\\
           \end{pmatrix},\;
         w_{\beta} = \begin{pmatrix}
            1  & 0 & 0  & 0\\
            0  & 0  & 0  & -1\\
            0  & 0  & 1  & 0\\
            0  & 1  & 0  & 0\\
           \end{pmatrix}.           
\end{align*}

We used results by {\sc Hua} and {\sc Curtis} \cite{Hua}, to extract an algorithm that expresses an arbitrary element $A \in Sp_4(\Z)$ as word in certain generators, which were then reexpressed into the {\sc Behr}-generators 
\[x_{\alpha},x_{\beta},x_{\alpha+\beta},x_{2\alpha+\beta},w_{\alpha},w_{\beta}.\]
For example the group the generators of $G(d,k)$ can be written as
\begin{align*}
 g_1 &= x_{\beta} \\
 g_2 &= (w_{\alpha} w_{\beta})^{-2} x_{2\alpha+\beta}^{-d}x_{\beta}^{k} x_{\alpha}^{-1}
	 w_{\alpha}^{-3}x_{\alpha}^{-1} (w_{\alpha} w_{\beta})^{-2}  
\end{align*}
Hence, if the generators of a finite index subgroup $M = \langle A_1,\ldots,A_n \rangle$ of $Sp_4(\Z)$ are given, we can try to use algorithms from computational group theory for finitely presented groups to compute the index $[Sp_4(\Z):M]$. In this way the results of theorem $3$ were found.\\ 

To understand Theorem $4$, one has to look a bit closer to the geometry 
associated to the finite symplectic group. It is a classical fact that 
$Sp_4(\Z/2)$, the reduction of $Sp_4(\Z)$ mod $2$,
is isomorphic to the permutation group $S_6$. A way to realise $Sp_4(\Z/2)$ naturally as a permutation group of six objects is the following. 
The $15$ points of $\P^3:=\P^3(\Z/2)$ correspond to the $15$ 
transpositions in $S_6$; the point pairs having symplectic scalar product equal to
one correspond to transpositions with a common index. The six five-tuples of 
transpositions all having a common index thus correspond to six five-tuples
of points in $\P^3$ that have pairwise symplectic scalar product equal to one.
Lets call such five-tuples a {\em pentade} of points.
These six pentades are permuted by $Sp_4(\Z/2)$, thus 
defining an isomorphism with the permutation group $S_6$. A subgroup
fixing such a pentade has index $6$ and is a copy of $S_5$.
Furthermore, there are $10$ {\em synthemes}, that is ways to divide 
six elements in two subsets of cardinality three. These correspond
however precisely to the pairs of triples of elements of $\P^3$ with the
property that the elements have symplectic scalar product one
if they belong to the same triple and zero else. The stabiliser of
such a syntheme is a subgroup of index $10$.

To make this explicit, let us label the elements of $\P^3$ by the letters
from a to o:

\[a=(0,0,0,1),\; b=(0,0,1,0),\; c=(0,0,1,1),\; d=(0,1,0,0),\]
\[e=(0,1,0,1),\; f=(0,1,1,0),\; g=(0,1,1,1),\; h=(1,0,0,0),\]
\[i=(1,0,0,1),\; j=(1,0,1,0),\; k=(1,0,1,1),\; l=(1,1,0,0),\]
\[m=(1,1,0,1),\; n=(1,1,1,0),\; o=(1,1,1,1) \]

One verifies at once that the six pentades are given by
\[ 1=\{a,d,g,m,o\},\; 2=\{a,e,f,l,n\},\; 3=\{b,h,k,n,o\},\]
\[ 4=\{b,i,j,l,m\},\; 5=\{c,d,e,i,k\},\; 6=\{c,f,g,h,j\}\]

These are permuted by $Sp_4(\Z/2)$. Indeed, a
transvection mod $2$ of an element $p \in \P^3$ 
\[ T_p: v \mapsto v+(v,p)p \]
acts as a transposition in the set $\{1,2,3,4,5,6\}$.
For example, one verifies that $T_a$ acts as the transposition
$(1,2)$. 
For the matrices with $d=k=1 \mod 2$ one finds
\[M \cdot a=a,\;M \cdot d=o,\;M \cdot g=m,\; M\cdot g=m,\; M\cdot m=d,\; M\cdot o=g\]
so that $M$ maps the pentade $1$ to it self,  $M \cdot 1=1$
In a similar way we obtain
\[M \cdot 1=1,\; M \cdot 2=2,\; M\cdot 3=6,\; M \cdot 4=5,\; M \cdot 5=3,\; M\cdot 6=4\]
\[N \cdot 1=5,\; N \cdot 2=2,\; N\cdot 3=3,\; N\cdot 4=4,\; N\cdot 5=1,\; N\cdot 6=6\]
so that only the pentade $2=\{a,e,f,l,n\}$ is fixed
by both $M$ and $N$ and one readily verifies that they generate
the stabiliser.

The ten synthemes, given as pairs of triples, are given by
\[I=\{\{a,d,e\},\{b,h,j\}\},\; II=\{\{a,f,g\},\{b,i,k\}\}\]
\[III=\{\{a,l,m\},\{c,h,k\}\},\; IV=\{\{a,n,o\},\{c,i,j\}\}\]
\[V=\{\{b,l,n\},\{c,d,g\}\},\; VI=\{\{b,m,o\},\{c,e,f\}\}\]
\[VII=\{\{d,i,m\},\{f,h,n\}\},\; VIII=\{\{d,k,o\},\{f,j,l\}\}\]
\[IX=\{\{e,i,l\},\{g,h,o\}\},\; X=\{\{e,k,n\},\{g,j,m\}\}\]
The group $Sp_4(\Z/2)$ permutes these synthemes, and one 
verifies that in case $d=1 \mod 2, k=0 \mod 2$ the matrix 
$M$ induces the permutation
\[ (I,IV,II,III)(VII,X,IX,IIIV)\]
and $N$ the permutation
\[(II,VI)(III,IX)(VI,X)\]
so that precisely syntheme $V=\{\{b,l,n\},\{c,d,g\}\}$ is preserved.\\

{\bf Remark:} There is another set of six objects that $Sp_4(\Z/2)$
permutes, which reflects the famous outer automorphism of $S_6$.
In the finite symplectic geometry these correspond to disjoint 
five-tuples of lagrangian lines. In the notation used above, these
are 
\[1'=\{\{a,b,c\},\{d,h,l\},\{e,j,o\},\{f,k,m\},\{g,i,n\}\},\]
\[2'=\{\{a,b,c\},\{d,j,n\},\{e,h,m\},\{f,i,o\},\{g,k,l\}\},\]
\[3'=\{\{a,j,k\},\{b,e,g\},\{c,m,n\},\{d,h,l\},\{f,i,o\}\},\]
\[4'=\{\{a,h,i\},\{b,d,f\},\{c,m,n\},\{e,j,o\},\{g,k,l\}\},\]
\[5'=\{\{a,h,i\},\{b,e,g\},\{c,l,o\},\{d,j,n\},\{f,k,m\}\},\]
\[6'=\{\{a,j,k\},\{b,d,f\},\{c,l,o\},\{e,h,m\},\{g,i,n\}\}\]

The stabiliser of such a pentade of lines is also isomorphic to $S_5$,
but is not conjugate to the stabiliser of a pentade of points.
The fact that the monodromy group $G(1,3)$ preserves a pentade of 
points rather than a pentade of lines is an intrinsic property and
is independent of any choices.

\section{An observation}

The dichotomy between cases of finite and infinite index
is rather mysterious. The finiteness of the index
does not seem to correlate to any simple geometrical
invariant of the Calabi-Yau. 
On the other hand, when we make the following plot the $14$ cases
in a diagram with where black boxes represent the cases of infinite index,
a pattern arises.
\begin{center}
\includegraphics[width=10cm]{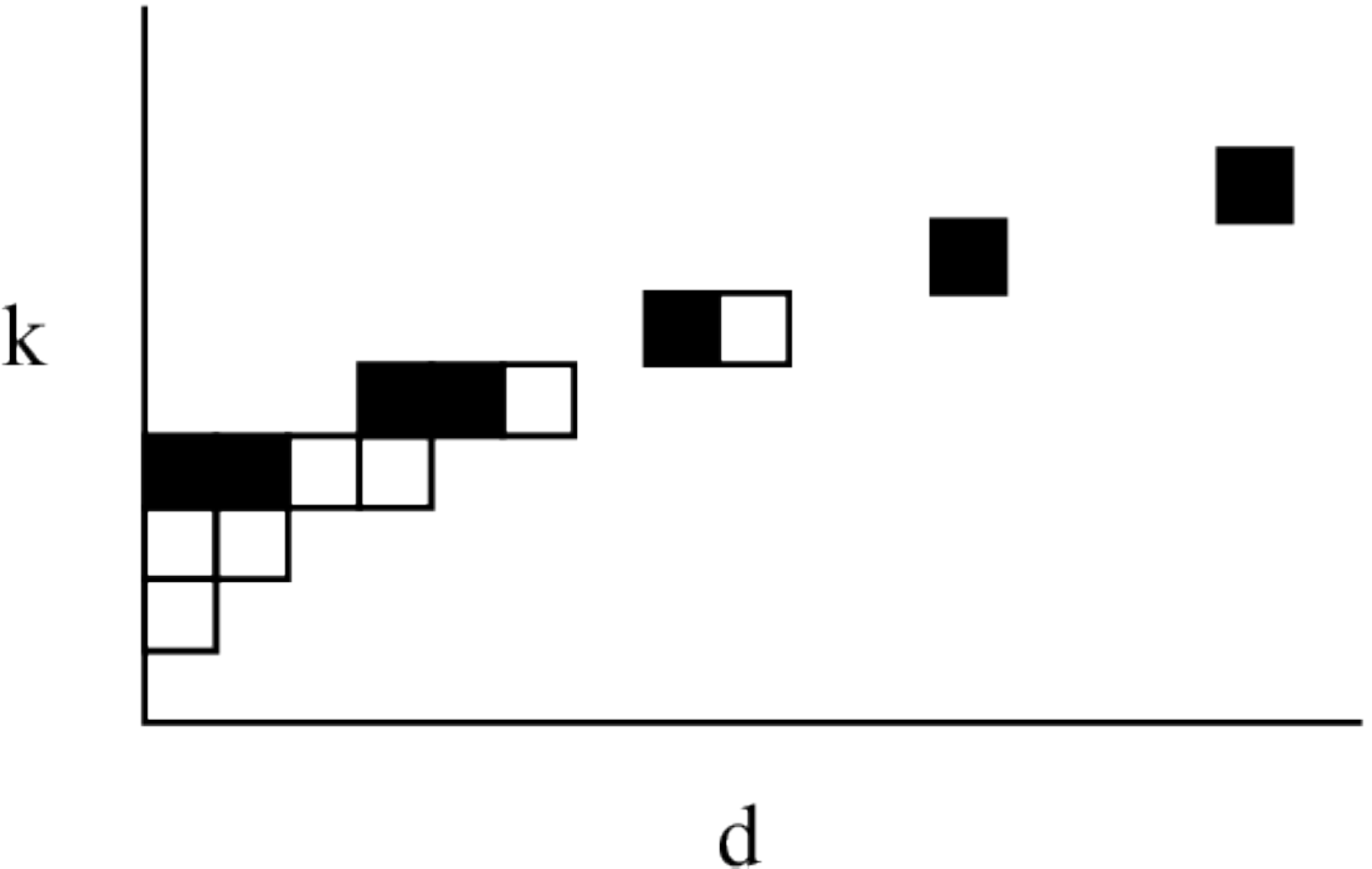}
\end{center}
There is a tendency for the finite index cases to be
lie ``under'' the infinite cases. Also, in the cases of finite index,
the index increases monotonously with $d$.
Apparently one may look at the quantity 
\[\Lambda:=\frac{7k-2d}{24}\]  
so that the cases with $\Lambda >1$ have infinite index and those
with $\Lambda <1$ have finite index. There are three cases where
$\Lambda=1$, to know $(2,4), (9,6), (16,8)$ of which only $(9,6)$
has finite index. 

\section{Non-Hypergeometric Operators with finite Index}

An obvious question is to ask in which cases of Calabi-Yau operators
from the list \cite{AESZ} have finite and which infinite index. Many of
these are ``conifold-operators'', which means that the singularity nearest
to the origin has exponents $0,1,1,2$. In such case one can define the 
invariants $d$ and $k$, and one is tempted to make the following\\

{\bf Wild Guess} Let $G \subset Sp_4(\Z)$ be the monodromy group of a conifold Calabi-Yau 
operator.
If $\Lambda >1$ then the index is infinite and if $\Lambda <1$ then the
index is finite.\\

Using this heuristic, we went through the list of Calabi-Yau operators
and discovered the following:\\

{\bf Theorem 5} The following {\em non-hypergeometric operators} have monodromy
of finite index in $Sp_4(\Z)$.
\[
\begin{array}{|c||c|c|c|c||c|c|}
\hline
\textup{AESZ}&H^3=d&k&c_2.H&c_3&\textup{Index}&G(d,k)-\textup{Index}\\
\hline
289&2&2&20&-16&360&5760\\
\hline
292&3&3&30&-92&6&933120\\
\hline
241&4&3&28&-60& 3840&122880\\
\hline
257&4&3&28&-32& 122880&122880\\
\hline
337&5&4&38&-102&1&3900000\\
\hline
33&6&4&24&-144&1036800&?\\
\hline
\end{array}
\] 
We included the index of the corresponding $G(d,k)$-group, as far as we could
determine it. Note that these groups do not belong to the family of $14$.
We note that the cases appearing here are all rather similar: all operators have apart from $0$ and $\infty$, two conifold points (exponents $0,1,1,2$) and a further apparent singularity (exponents $0,1,3,4$).\\
We list here the monodromy matrix around the extra conifold point in the basis explained in section 2. This monodomry transformation is also a symplectic reflection; we list the corresponding reflection vector.

\begin{longtable}{|c|c|c|}
\hline
\textup{Case}&\textup{Extra matrix}&\textup{Reflection vector}\\
\hline
289&$\begin{pmatrix}  -1& 4& 2& 2\\  -2& 5& 2& 2\\  -2& 4& 3& 2\\  4& -8& -4&-3\\ \end{pmatrix}$ & $(-2^{1/2},8^{1/2},2^{1/2},2^{1/2})$\\
\hline
292&$\begin{pmatrix}  0& 2& 1& 2\\  -2& 5& 2& 4\\  -1& 2& 2& 2\\  2& -4& -2&-3\\ \end{pmatrix}$&  $(-1,2,1,2)$ \\
\hline
241&$\begin{pmatrix}  -1& 2& 1& 2\\  -4& 5& 2& 4\\  -4& 4& 3& 4\\  4& -4& -2&-3\\ \end{pmatrix}$& $(-2,2,1,2)$\\
\hline
257&$\begin{pmatrix}  -3& 3& 1& 0\\  0& 1& 0& 0\\  -16& 12& 5& 0\\  12& -9&-3& 1\\ \end{pmatrix}$& $(-4,3,1,0)$\\
\hline
337&$\begin{pmatrix}  1& 0& 0& 0\\  1& 1& 0& 1\\ -1& 0& 1& -1\\  0& 0& 0&1\\  \end{pmatrix}$& $(1,0,0,1)$\\
\hline
33&$\begin{pmatrix}  1& 0& 0& 0\\  2& 1& 0& 2\\ -2& 0& 1& -2\\  0& 0& 0&1\\  \end{pmatrix}$ & $(2^{1/2},0,0,2^{1/2})$\\
\hline
\end{longtable}
Remarkable is the case $337$, which apparently has the full
$Sp_4(\Z)$ as monodromy group. The index of $G(5,4)$ is rather large, so in
this case the extra monodromy matrix makes a big difference. On the other hand,
for case $257$ the extra monodromy transformation does nothing, as in this
case the index is the same as for the group $G(4,3)$.\\
We believe that there are many more of cases of finite index in the list; this is currently under investigation. No geometrical incarnation of these operators
on the $A$-side is known to us, although we believe they should exist.\\

{\bf Operator AESZ 289 and Riemann-Symbol}\\
\begin{minipage}{9cm}\flushleft$\theta^4$
$-2^4x(400 \theta^4+2720 \theta^3+1752 \theta^2+392 \theta +33)$

$+2^{15}x^2(- 4272 \theta^4- 6288 \theta^3 + 3184 \theta^2+ 1484 \theta  +177)$

$+2^{24} 5 x^3(-4688 \theta^4  + 1536 \theta^3 + 1384 \theta^2+ 336 \theta +27)$ 
$+2^{36} 5^2 x^4(4\theta+1)(2\theta+1)^2(4\theta+3)$ 
\end{minipage} 
$\begin{Bmatrix}-\frac{1}{5120}&0&\frac{1}{16384}&\frac{1}{256}&\infty \\ \noalign{\smallskip} \hline
\begin{matrix} 0\\ \noalign{\smallskip}1\\ \noalign{\smallskip}3\\ \noalign{\smallskip}4 \end{matrix}&
\begin{matrix} 0\\ \noalign{\smallskip}0\\ \noalign{\smallskip}0\\ \noalign{\smallskip}0 \end{matrix}&
\begin{matrix} 0\\ \noalign{\smallskip}1\\ \noalign{\smallskip}1\\ \noalign{\smallskip}2 \end{matrix}&
\begin{matrix} 0\\ \noalign{\smallskip}1\\ \noalign{\smallskip}1\\ \noalign{\smallskip}2 \end{matrix}&
\begin{matrix}\frac{1}{4}\\ \noalign{\smallskip}\frac{1}{2}\\ \noalign{\smallskip}\frac{1}{2}\\ \noalign{\smallskip}\frac{3}{4}\end{matrix} \end{Bmatrix} $\\ 

{\bf Operator AESZ 292 and Riemann Symbol}\\
\begin{minipage}{9cm}\flushleft$9\theta^4$
$-2^2 3 x(4636 \theta^4+7928 \theta^3+ 5347 \theta^2 + 1383 \theta+  126)$
$+2^9 x^2(59048 \theta^4 + 50888 \theta^3 -26248\theta^2-16827 \theta-2205)$
$+2^{16} 7 x^3(- 9004 \theta^4+2304\theta^3+2511\theta^2+ 504 \theta +27)$
$-2^{24} 7^2 x^4 (4\theta+1)(2\theta+1)^2(4\theta+3)$ 
\end{minipage} 
$\begin{Bmatrix}-0.0853&0&0.000179&\frac{3}{896}&\infty \\ \noalign{\smallskip} \hline 
\begin{matrix} 0\\ \noalign{\smallskip}1\\ \noalign{\smallskip}1\\ \noalign{\smallskip}2 \end{matrix}&
\begin{matrix} 0\\ \noalign{\smallskip}0\\ \noalign{\smallskip}0\\ \noalign{\smallskip}0 \end{matrix}&
\begin{matrix} 0\\ \noalign{\smallskip}1\\ \noalign{\smallskip}1\\ \noalign{\smallskip}2 \end{matrix}&
\begin{matrix} 0\\ \noalign{\smallskip}1\\ \noalign{\smallskip}3\\ \noalign{\smallskip}4 \end{matrix}&
\begin{matrix}\frac{1}{4}\\ \noalign{\smallskip}\frac{1}{2}\\ \noalign{\smallskip}\frac{1}{2}\\ \noalign{\smallskip}\frac{3}{4}\end{matrix} \end{Bmatrix} $\\

{\bf Operator AESZ 241 and Riemann Symbol}\\
\begin{minipage}{9cm}\flushleft$\theta^4$
$-2^4x(152\theta^4+160\theta^3+110\theta^2+30\theta+3)$ 
$+2^{10} 3 x^2(428\theta^4+176\theta^3-299\theta^2-170\theta-25)$ 
$+2^{17} 3^2 x^3 (-136 \theta^4+216\theta^3+180\theta^2+51\theta+5)$ 
$-2^{24} 3^3 x^4 (3\theta+1)(2\theta+1)^2(3\theta+2)$ 
 \end{minipage} 
 $\begin{Bmatrix}-\frac{1}{64}&0&\frac{1}{1728}&\frac{1}{384}&\infty \\ \noalign{\smallskip} \hline
 \begin{matrix} 0\\ \noalign{\smallskip}1\\ \noalign{\smallskip}1\\ \noalign{\smallskip}2 \end{matrix}&
 \begin{matrix} 0\\ \noalign{\smallskip}0\\ \noalign{\smallskip}0\\ \noalign{\smallskip}0 \end{matrix}&
 \begin{matrix} 0\\ \noalign{\smallskip}1\\ \noalign{\smallskip}1\\ \noalign{\smallskip}2 \end{matrix}&
 \begin{matrix} 0\\ \noalign{\smallskip}1\\ \noalign{\smallskip}3\\ \noalign{\smallskip}4 \end{matrix}&
 \begin{matrix}\frac{1}{3}\\ \noalign{\smallskip}\frac{1}{2}\\ \noalign{\smallskip}\frac{1}{2}\\ \noalign{\smallskip}\frac{2}{3}\end{matrix} \end{Bmatrix} $
\\
{\bf Operator AESZ 257 and Riemann Symbol}\\
\begin{minipage}{9cm}\flushleft$\theta^4$
$-2^4 x (112\theta^4+416\theta^3+280\theta^2+72\theta+7)$ 
$+2^{12} x^2 (-656\theta^4-896\theta^3+216\theta^2+160\theta+23)$ 
$-2^{23} x^3 (96\theta^4+24\theta^3+12\theta^2+6\theta+1)$ 
$-2^{30} x^4 (2\theta+1)^4$ 
 \end{minipage} 
$\begin{Bmatrix}-0.0433&-\frac{1}{512}&0&0.000352&\infty \\ \noalign{\smallskip} \hline
 \begin{matrix} 0\\ \noalign{\smallskip}1\\ \noalign{\smallskip}1\\ \noalign{\smallskip}2 \end{matrix}&
 \begin{matrix} 0\\ \noalign{\smallskip}1\\ \noalign{\smallskip}3\\ \noalign{\smallskip}4 \end{matrix}&
 \begin{matrix} 0\\ \noalign{\smallskip}0\\ \noalign{\smallskip}0\\ \noalign{\smallskip}0 \end{matrix}&
 \begin{matrix} 0\\ \noalign{\smallskip}1\\ \noalign{\smallskip}1\\ \noalign{\smallskip}2 \end{matrix}&
 \begin{matrix}\frac{1}{2}\\ \noalign{\smallskip}\frac{1}{2}\\ \noalign{\smallskip}\frac{1}{2}\\ \noalign{\smallskip}\frac{1}{2}\end{matrix} \end{Bmatrix} $                                                                 

{\bf Operator AESZ 337 and Riemann Symbol}\\
\begin{minipage}{9cm}\flushleft$25\theta^4$
$-3 \cdot 5 x (3483 \theta^4 +6102 \theta^3+4241\theta^2+1190\theta+120)$ 
$+2^53^2 x^2 ( 31428 \theta^4+ 35559 \theta^3+ 243 \theta^2 - 4320 \theta-740 )$

$-2^8 3^5 x^3 (7371 \theta^4+4860 \theta^3+2997\theta^2+1080\theta+140)$
$+x^4 2^{13} 3^8 x^4 (3\theta+1)^2(3\theta+2)^2$ 
 \end{minipage} 
 $\begin{Bmatrix}0&0.000525&\frac{5}{432}&0.0816&\infty \\ \noalign{\smallskip} \hline 
 \begin{matrix} 0\\ \noalign{\smallskip}0\\ \noalign{\smallskip}0\\ \noalign{\smallskip}0 \end{matrix}&
 \begin{matrix} 0\\ \noalign{\smallskip}1\\ \noalign{\smallskip}1\\ \noalign{\smallskip}2 \end{matrix}&
 \begin{matrix} 0\\ \noalign{\smallskip}1\\ \noalign{\smallskip}3\\ \noalign{\smallskip}4 \end{matrix}&
 \begin{matrix} 0\\ \noalign{\smallskip}1\\ \noalign{\smallskip}1\\ \noalign{\smallskip}2 \end{matrix}&
 \begin{matrix}\frac{1}{3}\\ \noalign{\smallskip}\frac{1}{3}\\ \noalign{\smallskip}\frac{2}{3}\\ \noalign{\smallskip}\frac{2}{3}\end{matrix} \end{Bmatrix} $                                                                 

{\bf Operator AESZ 33 and Riemann Symbol}\\
\begin{minipage}{9cm}\flushleft$\theta^4$
$-2^2 x (324\theta^4+456\theta^3+321\theta^2+93\theta+10)$ 
$+2^9 x^2 (584 \theta^4+584\theta^3+4\theta^2-71\theta-13)$ 
$-2^{16} x^3 (324 \theta^4+192 \theta^3+123\theta^2+48\theta+7)$ 
$+2^{24} x^4 (2\theta+1)^4$ 
 \end{minipage} 
 $\begin{Bmatrix}0&\frac{1}{1024}&\frac{1}{128}&\frac{1}{16}&\infty \\ \noalign{\smallskip} \hline 
 \begin{matrix} 0\\ \noalign{\smallskip}0\\ \noalign{\smallskip}0\\ \noalign{\smallskip}0 \end{matrix}&
 \begin{matrix} 0\\ \noalign{\smallskip}1\\ \noalign{\smallskip}1\\ \noalign{\smallskip}2 \end{matrix}&
 \begin{matrix} 0\\ \noalign{\smallskip}1\\ \noalign{\smallskip}3\\ \noalign{\smallskip}4 \end{matrix}&
 \begin{matrix} 0\\ \noalign{\smallskip}1\\ \noalign{\smallskip}1\\ \noalign{\smallskip}2 \end{matrix}&
 \begin{matrix}\frac{1}{2}\\ \noalign{\smallskip}\frac{1}{2}\\ \noalign{\smallskip}\frac{1}{2}\\ \noalign{\smallskip}\frac{1}{2}\end{matrix} \end{Bmatrix} $   

\section{Monodromy group mod $N$}
Using {\tt GAP}, we can also try to determine the structure of the monodromy
group in $Sp_4(\Z/N\Z)$ for various $N$. Note
that 
\[ |Sp_4(\Z/N\Z)|=N^{10}\prod_{p \mid N }(1-p^{-2})(1-p^{-4})\]
For convenience of the reader we list the result of a {\tt GAP}-computation.

\setlength{\tabcolsep}{0.5mm}
{\tiny
\[
\begin{tabular}{|r|r|r|r|r|r|r|r|r|r|r|r|r|r|r|r|}
\hline
N&(1,4)&(1,3)&(1,2)&(2,4)&(2,3)&(3,4)&(4,5)&(4,4)&(5,5)&(6,5)&(8,6)&(9,6)&(12,7)&(16,8)\\
\hline
\hline
 2& 10&6&10&   90& 60&    10&   60&     90&    6&     60&     90&    10&      60&      90\\
 3&  1&1& 1&    1&  1&   720&    1&      1&    1&    720&      1&   640&     720&       1\\
 4&160&6&10& 2880&240&   160& 3840&   5760&    6&    240&   5760&    10&    3840&    5760\\
 5&  1&1& 1&    1&  1&     1&    1&      1&14976&      1&      1&     1&       1&       1\\
 6& 10&6&10&   90& 60&  7200&   60&     90&    6&  43200&     90&  6400&   43200&      90\\
 7&  1&1& 1&    1&  1&     1&    1&      1&    1&      1&      1&     1&       1&       1\\
 8&160&6&10&46080&960&   160&15360& 184320&    6&    960& 368640&     10&  15360&  368640\\
 9&  1&1& 1&    1&  1& 19440&    1&      1&    1&  19440&      1& 466560&  19440&       1\\
10& 10&6&10&   90& 60&    10&   60&     90&89856&     60&     90&     10&     60&      90\\
11&  1&1& 1&    1&  1&     1&    1&      1&    1&      1&      1&      1&      1&       1\\
12&160&6&10& 2880&240&115200& 3840&   5760&    6& 172800&   5760&   6400&2764800&    5760\\
13&  1&1& 1&    1&  1&     1&    1&      1&    1&      1&      1&      1&      1&       1\\
14& 10&6&10&   90& 60&    10&   60&     90&    6&     60&     90&     10&     60&      90\\
15&  1&1& 1&    1&  1&   720&    1&      1&14976&    720&      1&      1&    640&     720\\
16&160&6&10&92160&960&   160&61440&2949120&    6&    960&5898240&     10&  61440&23592960\\
17&  1&1& 1&    1&  1&     1&    1&      1&    1&      1&      1&      1&      1&       1\\
18& 10&6&10&   90& 60&194400&   60&     90&    6&1166400&     90&4665600&1166400&      90\\
19&  1&1& 1&    1&  1&     1&    1&      1&    1&      1&      1&      1&      1&       1\\
20&160&6&10& 2880&240&   160& 3840&   5760&89856&    240&   5760&     10&   3840&    5760\\
21&  1&1& 1&    1&  1&   720&    1&      1&    1&    720&      1&    640&    720&       1\\
22& 10&6&10&   90& 60&    10&   60&     90&    6&     60&     90&     10&     60&      90\\
23&  1&1& 1&    1&  1&     1&    1&      1&    1&      1&      1&      1&      1&       1\\
24&160&6&10&46080&960&115200&15360& 184320&    6& 691200& 368640&   6400&11059200& 368640\\    
25&  1&1& 1&    1&  1&     1&    1&      1&46800000&   1&      1&      1&      1&       1\\
26& 10&6&10&   90& 60&    10&   60&     90&    6&     60&     90&     10&     60&      90\\
27&  1&1& 1&    1&  1& 19440&    1&      1&    1&  19440&      1&113374080&19440&       1\\
\hline
\end{tabular}
\]
}

The table contains some redundancies: if $N$ and $M$ have no common
factor, the index in $Sp_4(\Z/NM)$ is the product of the indices 
in $Sp_4(\Z/N)$ and $Sp_4(\Z/M)$. The table also shows some remarkable 
phenomena. The case $(1,4)$ is of infinite
index in $Sp_4(\Z)$, but the reductions mod $N$ suggest the index is $160$
when considered $2$-adically, that is in the group $Sp_4(\Z_2)$. The columns $(1,3)$, $(1,2)$, $(2,3)$ look very similar, but here the index in $Sp_2(\Z)$ indeed is $6$, $10$, $960$, respectively. For $(5,5)$ the numbers probably will grow further;
note the prime number $13$ entering in the index. All other columns have
only $2$, $3$ and $5$ appearing in the prime factorisation. The column $(9,6)$
shows that the index in the last case of finite index is at least 
\[ 90 \cdot 113374080=10203667200=2^83^{13}5^2\]  
and might very well be equal to this number.\\

{\bf Acknowledgement:}
This work was begun in october $2012$ during a stay of the authors at the MSRI
in Berkeley. We thank this institution for its hospitality and its excellent 
research environment. We thank S. Singh, H. Thomas and W. Zudilin for showing
interest in early versions of this work. Thanks also to C. Doran, who has
indicated that he has work in preparation that is related to this paper.

\end{document}